\documentstyle[11pt,amssymb,amstex]{amsart}

\newtheorem{thm}{Theorem}[section]
\newtheorem{lem}[thm]{Lemma}
\newtheorem{cor}[thm]{Corollary}

\begin{document}

\title[On dissipative quadratic stochastic operators]
{On dissipative quadratic stochastic operators}

%\thanks{On leave from Department of Mechanics and Mathematics,
%National University of Uzbekistan, Tashkent, 100174, Uzbekistan}

\author{Farruh Shahidi}
\address{Farruh Shahidi \\
Department of Mechanics and Mathematics,\\
National University of Uzbekistan\\
Vuzgorodok, Tashkent, 700174, Uzbekistan} \email{{\tt
farruh.shahidi@@gmail.com}}

\begin{abstract}
In present paper we introduce the notion of dissipative quadratic
stochastic operator and cubic stochastic operator. We prove
necessary conditions for dissipativity of quadratic stochastic
operators. Besides, it is studied certain limit behavior of such
operators. Finally we prove ergodic  theorem for dissipative
operators.\vskip 0.3cm \noindent {\it
Mathematics Subject Classification}: 15A51, 47H60, 46T05, 92B99.
\end{abstract}

\maketitle

\section{Introduction}

It is known \cite{Lyu} that the theory of quadratic stochastic
operators frequently arise in many models of physics, biology and
so on. Let us shortly mention how such kind of operators appear in
population genetics. Consider biological population, that is, a
community of organisms closed with respect to reproduction
\cite{Ber}.Assume that every individual in this population belongs
to one of the species $1,2,\cdots, m.$ The scale of species in
such that the species of parents $i$ and $j$ unambiguously
determine the probability of every species $k$ for the first
generation of direct descendants. We denote this probability (the
heredity coefficient) via $p_{ij,k}$ and
$\sum\limits_{k=1}^{m}p_{ij,k}=1$ for all $i,j.$ Assume that the
population so large that frequency fluctuations can be neglected.
Then the state of the population can be described by the tuple
$x=(x_{1},x_{2},\cdots,x_{m})$ of species probabilities, that is,
$x_{i}$ is the fraction of the species $i$ in the population. In
the case of panmixia(random interbreeding), the parent pairs $i$
and $j$ arise for a fixed state $x=(x_{1},x_{2},\cdots,x_{m})$
with probability $x_{i}x_{j}.$ Hence
$$x_{k}'=\sum\limits_{i,j=1}^{m}p_{ij,k}x_{i}x_{j}$$ is the total
probability of the species in the first generation of direct
descendants.  Note that the concept of quadratic stochastic
operator firstly introduced by Bernstein in \cite{Ber}. A lot of
papers were devoted to study such operators (see for example,
\cite{Ga1},\cite{Lyu},\cite{K},\cite{Val},\cite{M3},\cite{SG}).
One of the central problem in this theory is the study limit
behavior and ergodic properties of trajectories of quadratic
operators (\cite{GZ},\cite{M1},\cite{M2},\cite{SG}).Note that the
studying of a such properties of quadratic stochastic operators is
very difficult. Even in the two-dimensional simplex the problem is
still unsolved.This problem is well studied for Volterra quadratic
stochastic operators(\cite{Ga1}). We show that operators, which we
will study is not Volterra operators. In \cite{Ga2} a class of so
called bistochastic operators, i.e. operators with property
$Vx\prec x$ for all $x\in S^{m-1}$, is outlined (see next section
for notations).It should be mentioned that bistochastic operators
do not dissipate, i.e. no species of such operators presented in
the beginning of the evolution can be disappear. In the present
paper we are going to study quadratic and cubic operators
satisfying the condition $Vx\succ x$ for every $x\in S^{m-1}$,
which will be called {\it dissipative}. It should be mentioned
that the class of bistochastic and dissipative operators are
almost disjoint, i.e. only permutation operators can be
bistochastic as well as dissipative.

The paper organized as follows. In section 2 we give some
preliminaries on quadratic stochastic operators and definition of
dissipative ones. There we show that the set of dissipative
operators does not form a convex set, while the class of
bistochastic operators is convex. In section 3 we study certain
limit properties of dissipative ones. Moreover, we describe all
such operators in small dimensions. In section 4 we prove that
every dissipative operator satisfies an ergodic theorem. In
section 5 we study dissipative cubic stochastic operators.Finally,
in section 6 we give a conclusion of obtained results.

\section{Preliminaries}

Let $S^{m-1}=\{x\in R^{m}: x_{i}\geq 0, \ \
\sum\limits_{i=1}^{m}x_{i}=1\}$ be a $(m-1)-$dimensional simplex.
Let $e_k=(0,0,\cdots,\underbrace{1}_{k},\cdots, 0)$,($
k=\overline{1,m}$) be its vertices.

An operator $V:S^{m-1}\rightarrow S^{m-1}$ is called a stochastic
operator. A stochastic operator is called a {\it quadratic
stochastic operator (q.s.o. in short)} if

$$(Vx)_{k}=\sum\limits_{i,j=1}^{m}p_{ij,k}x_{i}x_{j}.\eqno (1)$$

Here the coefficients $p_{ij,k}$ satisfy the following conditions
$$p_{ij,k}=p_{ji,k}\geq 0 , \ \ \sum\limits_{k=1}^mp_{ij,k}=1.\eqno (2)$$
It is easy to see that q.s.o. is well defined, i.e. it maps the
simplex into itself.

An operator $V:S^{m-1}\rightarrow S^{m-1}$ is called {\it cubic
stochastic operator(c.s.o. in short)} if it has a following form
$$(Vx)_k=\sum\limits_{i,j,l=1}^{m}p_{ijl,k}x_ix_jx_l \ \ k=\overline {1,m},\eqno (3)$$
where       $x=(x_1,x_2,...,x_m)\in S^{m-1},$ and
$$ p_{ijl ,k}=p_{ilj,k}=p_{jil,k}=p_{jli,k}=p_{lij,k}=p_{lji ,k} \geqslant 0, \ \forall i,j,l,k=\overline {1,m}$$
$$ \sum\limits_{k=1}^{m}p_{ijl,k}=1\ \forall
i,j,l=\overline {1,m}  .$$

Now, let $x,y\in R^m.$ Let's put
$x_{\downarrow}=(x_{[1]},x_{[2]},\cdots x_{[m]}),$ where
$(x_{[1]},x_{[2]},\cdots x_{[m]})$- decreasing rearrangement of
$(x_1,x_2,\cdots x_m),$ that is $x_{[1]}\geq x_{[2]}\geq \cdots
\geq x_{[m]}.$

\textbf{Definition 1.} We say that $x$ {\it majorized} by $y$ (or
$y$ {\it majorates} $x$), and write $x\prec y$(or $y\succ x$) if
the following conditions are fulfilled:

1) $ \sum\limits_{i=1}^{k}x_{[i]}\leq
\sum\limits_{i=1}^{k}y_{[i]}, \ \ k=\overline{1,m-1}$

2) $ \sum\limits_{i=1}^{m}x_{[i]}= \sum\limits_{i=1}^{m}y_{[i]}$

\begin{lem}\label{1} \cite{Ma}  For any
$x=(x_1,x_2,\cdots, x_m)\in S^{m-1}$ we have
$$(\frac{1}{m},\frac{1}{m},\cdots ,\frac{1}{m})\prec x\prec (1,0,\cdots 0).$$
\end{lem}

\textbf{Remark.}  It should be noted that  $"\prec"$ is not a
partial ordering, because from $x\prec y$ and $y\prec x$ it only
follows that  $x_{\downarrow}=y_{\downarrow}.$

\textbf{Definition 2.} A stochastic operator  $V$ is called {\it
dissipative} if
$$
Vx\succ x, \ \ \forall x\in S^{m-1}.\eqno (4)
$$
In particular q.s.o.(1) is called dissipative q.s.o. if it
satisfies the above mentioned condition (4). By this analogue we
define dissipative cubic stochastic operators.

\textbf{Observation} Consider the case when $V$ is linear
dissipative operator, that is $Vx=Ax$, here
$A=(a_{ij})_{i,j=\overline{1,m}}.$ Now we show that only
permutation linear operators are dissipative.

Since $Vx\succ x$ then by putting $x=e_i$ we have $Ae_i \succ
e_i.$ From the lemma 2.1. we obtain that
$(Ae_i)_{\downarrow}=(e_i)_{\downarrow}.$ The last is mean that
only one component of the vector $Ae_i$ is $1$ an the others are
$0$, which means that dissipative linear operators are permutation
operators. Therefore in linear case studying dissipative operators
is very simple. More interesting to study non-linear dissipative
operators.

We mention that dissipative operators are not studied, except so
called F- quadratic stochastic operators \cite{Roz}. Let us
consider the set $E=\{0,1,\cdots, m-1\}.$ Fix a set $F\subset E$
and call the set of "females" and the set $M=E\setminus F$ the set
of "males". The element $0$ will play the role of empty body
$$
p_{ij,k}=\left\{
\begin{array}{ll}
1,   \ \ \ \ if \  k=0, i,j\in F\cup {0} \  or  \ \ i,j\in M\cup {0}\\
0,   \ \ \ \ if \  k\neq 0, i,j\in F\cup {0} \  or  \ \ i,j\in M\cup {0}\\
\ge 0, \ \ if \ i\in F, j\in m, \forall k.
\end{array}
\right.
$$

Biological treatment of the above coefficients is very clear: a
"child" k can be generated if its parents are taken from different
classes $F$ and $M.$ In general $p_{ij,0}$ can be strictly
positive for $i\in F$ and $j\in M$, this corresponds, for
instance, to the case when "female" $i$ with "male" $j$ can not
generate a "child" since one of them(or both) is(are) ill. In
general, if we set $F=\{1\}$ then  F- q.s.o. have the following
view:

$$V:\left\{
\begin{array}{ll}
V(x_0)=1-2x_1\sum\limits_{i=2}^m(1-p_{1i,0})x_i,\\
V(x_k)=2x_1\sum\limits_{i=2}^mp_{1i,k}x_i, \ \ k=1,2,\cdots, m-1.
\end{array}
\right.$$

In \cite{Roz}, authors studied limit behavior of such operators
for arbitrary set $F$. There was proven that such operator has a
unique fixed and trajectory tend to this point. Only for the case
of $F=\{1\}$ our class of operators and F- q.s.o. can intersect.
In the next section we study limit behavior of dissipative ones,
and show that our class of operators is wider in some sense.

Now recall the term of well known Volterra q.s.o. A q.s.o. (1) is
called Volterra q.s.o. if it satisfies an additional assumption
$p_{ij,k}=0 \forall k\notin \{i,j\}$.  By changing
$a_{ki}=2p_{ik,k}-1$ one can write down the following canonical
form:
$$(Vx)_k=x_k(1+\sum\limits_{i=1}^ma_{ki}x_i)$$
In \cite{Ga1}, it was proved that for any non-fixed initial point
from the interior of the simplex, trajectory  approaches a bound
of the simplex. Also we show that dissipative q.s.o. can not be
Volterra q.s.o.

Now let us introduce the last notation which will be useful for
the next sections.

The point $x^0\in S^{m-1}$ is called fixed if $Vx^0=x^0.$ As a
rule there are tree types of fixed point.

We call a fixed point $x^0$ elliptic (hyperbolic; parabolic) if
the spectrum of Jacobian $J(x^0)$ restricted to the invariant
plane {$\sum\limits_{i=1}^{m}x_i=0 $} lies inside the unit ball
(respectively, outside the closure of the unit ball; inside the
unit circle).

\section{ Dissipative quadratic stochastic operators and their limit behavior}

In this section we are going to study the regularity of certain
dissipative q.s.o. Also we give an example of dissipative q.s.o.
which has infinitely many fixed points and study its limit
behavior.

The following example shows that the set of dissipative q.s.o. is
non-empty.
\begin{eqnarray*}
&&(Vx)_1=x_1^2+x_2^2+x_3^2+x_1x_2+x_1x_3+x_2x_3,\\
&&(Vx)_2=x_1x_2+x_1x_3,\\
&&(Vx)_3=x_2x_3.
\end{eqnarray*}

Let us  denote $a_{ij}=(p_{ij,1},p_{ij,2},\cdots p_{ij,m})\ \
\forall i,j=\overline {1,m}$ where $p_{ij,k}$ are the coefficients
of q.s.o. (1). One can see that $a_{ij}\in S^{m-1}$, for all
$i,j\in\overline{1,m}$

\begin{lem}\label{2} Let $V$ be a dissipative q.s.o.
Then the following conditions hold
$$ (a_{ii})_{\downarrow} =e_1 \ \ \forall i=\overline{1,m}.$$
\end{lem}

\begin{pf} We have $Vx\succ x \ \ \forall x\in S^{m-1}.$
By putting $x=e_i$ we get $e_i\prec Ve_i.$ On the other hand from
Lemma \ref{1} it follows that $e_i\succ x \ \ \forall x\in
S^{m-1}.$ That's why $(e_i)_{\downarrow}=(Ve_i)_{\downarrow}.$
Then the equality $Ve_i=a_{ii}$ implies the assertion.
\end{pf}

{\bf Remark.} Note that in \cite{Ga2} quadratic bistochastic
operators were studied, that is, operators satisfying the
condition $x\succ Vx \ \ \forall x\in S^{m-1}$. It was proved that
such operators form a convex compact set and its extreme points
were studied. The situation under consideration is different.
Indeed, let us consider the following operators:
\begin{eqnarray*}
&&(V_0x)_1= x_1x_2+x_1x_3,\\
&&(V_0x)_2=x_1^2+x_2^2+x_3^2+x_1x_2+x_2x_3+x_1x_3,\\
&& (V_0x)_3= x_2x_3.
\end{eqnarray*}
\begin{eqnarray*}
&&(V_1x)_1=x_1^2+x_2^2+x_3^2+x_1x_2+x_2x_3+x_1x_3,\\
&&(V_1x)_2=x_1x_2+x_1x_3,\\
&&(V_1x)_3=x_2x_3
\end{eqnarray*}
One can see that these operators are dissipative. However, Lemma
\ref{2} implies that the operator $V_{\lambda}=\lambda
V_1+(1-\lambda)V_0$ is not dissipative for any $\lambda\in(0,1)$.
Hence, the set of all dissipative q.s.o. is not convex.

Let $V$ be a dissipative q.s.o. Then thanks to Lemma \ref{2} it
can represented by
$$(Vx)_k=\sum\limits_{i\in
\alpha_k}x_i^2+2\sum\limits_{i<j}p_{ij,k}x_ix_j \ \
k=\overline{1,m},\eqno (5)$$ where $\alpha_k\subset
I=\{1,2\cdots,m\},$ $\alpha_i\cap \alpha_j=\emptyset,$ $i\neq j,$
$\bigcup\limits_{k=1}^m\alpha_k=I.$

\begin{lem}\label{3} Let (5) be a dissipative q.s.o.
\begin{itemize}
\item[(i)] If $j\in \alpha_{k_0},$ then
$p_{ij,k_0}=(a_{ij})_{[1]}\geq \frac{1}{2}$ $\forall
i=\overline{1,m}.$ \item[(ii)] If $m\geq 3,$ then
$(a_{ij})_{[k]}=0$ $\forall k\geq 3,$ $\forall i=\overline{1,m}.$
\end{itemize}
\end{lem}

\begin{pf}  (i).  Let $j\in \alpha_{k_0}$ and
$x=(1-\lambda)e_j+\lambda e_i.$ Here, as before, $e_i, e_j$ are
the vertices of the simplex and $\lambda$ is sufficiently small
positive number. It is easy to see that $x_{[1]}=1-\lambda$ and
$(Vx)_{[1]}=(Vx)_{k_0}.$ Since $Vx\succ x$ then $x_{[1]}\leq
(Vx)_{[1]},$ so $1-\lambda\leq (Vx)_{k_0}$ or
$$1-\lambda\leq (1-\lambda)^2+2p_{ij,k_0}\lambda(1-\lambda).$$
The last inequality implies that $p_{ij,k_{0}}\geq \frac{1}{2}.$

(ii). Denote $p_{ij,k^\ast}=\max\limits_{t\neq k_{0}}p_{ij,t}.$
One can see that $(Vx)_{k^{\ast}}=(a_{ij})_{[2]}.$ Now from
$$x_{[1]}+x_{[2]}\leq Vx_{[1]}+Vx_{[2]}.$$
 we obtain
 $$1\leq
(1-\lambda)^2+2(p_{ij,k_o}+p_{ij,k^{\ast}})\lambda(1-\lambda).$$
Assume that $\lambda\rightarrow 0.$ Then one gets
$p_{ij,k_o}+p_{ij,k^{\ast}}\geq 1.$ This yields that
$p_{ij,k_o}+p_{ij,k^{\ast}}=1$ and $(a_{ij})_{[k]}=0$ $\forall
k\geq 3,$ $\forall i=\overline{1,m}.$
\end{pf}

{\bf Observation. } Note that statements of Lemmas \ref{2} and
\ref{3} are the necessary conditions for q.s.o. to be dissipative.
It turns out that at $m=2$ the statements are sufficient. Indeed,
in this case only dissipative q.s.o. are  the identity operator
and the following one
\begin{eqnarray*}
&&(Vx)_1=x_1^2+x_2^2+ax_1x_2,\\
&&(Vx)_2=(2-a)x_1x_2, \end{eqnarray*} up to permutation of the
coordinates. Here $1\leq a\leq 2.$

However, when $m\geq 3$ then the statements are not sufficient.
Consider the following example of q.s.o.
\begin{eqnarray*}
&&(Vx)_1=x_1+x_2-x_1x_2,\\
&&(Vx)_2=0.8x_1x_2,\\
&&(Vx)_3=x_3+0.2x_1x_2.
\end{eqnarray*}
One can see that it satisfies the mentioned statements. But for
$x^0=(0.5;0.49;0.01)$ we have $Vx^0\nsucc x^0$, which means that
it is not dissipative.

Studying of limit behavior of all dissipative q.s.o. is a
difficult problem. We consider some particular cases. First recall
a q.s.o. $V:S^{m-1}\rightarrow S^{m-1}$ is called {\it regular} if
the trajectory of any $x\in S^{m-1}$ converges to a unique fixed
point.

Note that regular operators a'priori must have unique fixed point
and its fixed point is attracting.

Now let consider the case $\alpha_1=I$ and $\alpha_k=\emptyset$
for $k\neq 1.$ Then the operator has the following form
$$
\left.
\begin{array}{ll}
(Vx)_1=\sum\limits_{i=1}^mx_i^2+2\sum\limits_{i<j}p_{ij,1}x_ix_j\\
(Vx)_k=2\sum\limits_{i<j}p_{ij,k}x_ix_j, \ 2\leq k\leq m\\
\end{array}
\right\} \ \  \eqno (6)
$$

\begin{thm}\label{4}  A q.s.o. given by (6)  is regular. Its
unique fixed point is $e_1.$
\end{thm}

\begin{pf} Let us first prove that there is a unique
fixed point. The existence one follows from the Bohl-Brower
theorem. Denote it by
$x^{(0)}=(x^{(0)}_1,x^{(0)}_2,\cdots,x^{(0)}_m).$ It is clear that
tt satisfies the following equality:
$$x^{(0)}_1=\sum\limits_{i=1}^m(x^{(0)}_i)^2+2\sum\limits_{i<j}p_{ij,1}x^{(0)}_ix^{(0)}_j,$$
which can be rewritten by
$$(x^{(0)}_1)^2+x^{(0)}_1x^{(0)}_2+\cdots+x^{(0)}_1x^{(0)}_m=\sum\limits_{i=1}^m(x^{(0)}_i)^2+
2\sum\limits_{i<j}p_{ij,1}x^{(0)}_ix^{(0)}_j$$ or
$$(x^{(0)}_2)^2+(x^{(0)}_3)^2+\cdots+(x^{(0)}_m)^2+\sum\limits_{i=2}^m(2p_{ij,1}-1)x^{(0)}_1x^{(0)}_i+
\sum\limits_{2\leq i<j}p_{ij,1}x^{(0)}_ix^{(0)}_j=0.$$ It follows
from Lemma \ref{3} that $2p_{ij,1}-1\geq 0.$ Therefore, the left
hand side is positive, which means that the equality holds iff
$x_2=x_3=\cdots=x_m=0.$ Hence, q.s.o. (6) has a unique fixed point
$e_1.$

Now let us show that the operator is regular. Consider a function
$\varphi:S^{m-1}\rightarrow R,$ defined by
$$
\varphi(x)=x_2+x_3+\cdots+x_m.
$$
Then
$$\varphi(Vx)=\sum\limits_{i<j}\sum\limits_{k=2}^{m}2p_{ij,k}x_ix_j.$$

One can see that $2\sum\limits_{k=2}^mp_{ij,k}\leq 1$, since
$2p_{ij,1}-1\geq 0.$  Hence,
$$\varphi(Vx)\leq \sum\limits_{i<j}x_ix_j\leq
\sum\limits_{i=2}^mx_i\sum\limits_{i=1}^mx_i=\sum\limits_{i=2}^mx_i=\varphi(x).$$
Consequently, $\{\varphi(V^nx)\}$ is a decreasing sequence.
Therefore it converges.

Denote
$$\lim\limits_{n\rightarrow
\infty}\varphi(V^nx)=C.$$

The equality
$$(V^{n+1}x)_1=\sum\limits_{i=1}^m((V^nx)_i^2+2\sum\limits_{i<j}p_{ij,1}(V^nx)_i^{(n)}(V^nx)_j$$
with $2p_{ij,1}-1\geq 0$ implies
\begin{eqnarray*}
(V^{n+1}x)_1&=&\sum\limits_{i=1}^m((V^nx)_i^2+2\sum\limits_{i<j}p_{ij,1}(V^nx)_i(V^nx)_j\\
&\geq&
\sum\limits_{i=1}^m((V^nx)_i)^2+\sum\limits_{i<j}(V^nx)_i(V^nx)_j\\
&=&(V^nx)_1+\sum\limits_{i=2}^m((V^nx)_i)^2+
\sum\limits_{1<i<j}^m(V^nx)_i(V^nx)_j.
\end{eqnarray*}

Since $\lim\limits_{n\rightarrow
\infty}(V^{n+1}x)_1=\lim\limits_{n\rightarrow
 \infty}(V^nx)_1=1-C$ then we have
 $$\lim\limits_{n\rightarrow\infty}(\sum\limits_{i=2}^m((V^nx)_i)^2+
\sum\limits_{1<i<j}^m(V^nx)_i(V^nx)_2\leq 0.$$ The last means that
$\lim\limits_{n\rightarrow \infty}(V^nx)_i=0$ for all $i\geq 2,$
hence $\lim\limits_{n\rightarrow \infty}(V^nx)=e_1.$
 Therefore we deduce that the operator (6) is regular.
\end{pf}

\textbf{Remark} If the trajectory of q.s.o. belongs to the edge of
the simplex, then it  means that in process of time some species
of the population in the bound of disappearing.In our case we can
conclude that almost all species will disappear.

\textbf{Remark.} Note that q.s.o. of the form (6) need not to be
dissipative. Indeed, consider the following example.
\begin{eqnarray*}
&&(V_1x)_1=x_1^2+x_2^2+x_3^2+a_1x_1x_2+b_1x_2x_3+c_1x_1x_3,\\
&&(V_1x)_2=a_2x_1x_2+b_2x_2x_3+c_2x_1x_3,\\
&&(V_1x)_3=a_3x_1x_2+b_3x_2x_3+c_3x_1x_3
\end{eqnarray*}
here $a_1,b_1,c_1\geq \frac12$ and
$\sum\limits_{i=1}^3a_i=\sum\limits_{i=1}^3b_i=\sum\limits_{i=1}^3c_i=2$
If the coefficient are strictly positive, then q.s.o. need not to
be dissipative. Therefore, Theorem 3.1 is valid not only for
dissipative ones. On the other hand, we mentioned that in
\cite{Roz} a F-q.s.o. has been studied. Such operators can be
represented in (6) form, for the case of $F=\{1\}$. So, our
Theorem 3.1 is a generalization of a result in \cite{Roz} for the
case of $F=\{1\}$.

\begin{cor} In the case of $\alpha_k=I$ and
$\alpha_j=\emptyset$ for $j\neq
 k$ dissipative q.s.o. is regular and has a unique fixed point $e_k.$
 \end{cor}

Let $V$ be a q.s.o. Then  the set $\omega(x^0)=\bigcap \limits_{k
\geq 0}\overline{\bigcup \limits_{n\geq k} \bigl \{ V^nx^0 \bigr
\}}$ is called {\it $\omega$-limit set} of trajectory of initial
point $x^0\in S^{m-1}$. From the compactness of the simplex one
can deduce that $\omega(x^0)\neq \emptyset$ for all $x^0\in
S^{m-1}$.

Now we turn to the another case, namely let
$\alpha_1=I\backslash\{l\}$ $\alpha_2=\{l\}$ (actually we can put
$\alpha_{k_{0}}=\{l\}$ for some $k_{0}$) and $\alpha_k=\emptyset \
\ \forall k\geq 3.$ Then operator (5) has the following form
$$
\left.
\begin{array}{lll}
(Vx)_1=\sum\limits_{i=1, \ i\neq
l}^mx_i^2+2\sum\limits_{i<j}p_{ij,1}x_ix_j\\[5mm]
(Vx)_2=x_l^2+2\sum\limits_{i<j}p_{ij,2}x_ix_j\\[5mm]
(Vx)_k=2\sum\limits_{i<j}p_{ij,k}x_ix_j, \ \ 3\leq k\leq m
\end{array}
\right\} \eqno(7)
$$

\begin{thm}\label{5} If $l\neq 2$ then the operator (7) is
regular and has a unique fixed point $e_1$. If $l=2$ then the
operator (7) has infinitely many fixed points and all of them are
parabolic. Moreover, $\omega$-limit set of trajectory of any
initial point $x^0$ belongs to $co\{e_1,e_2\},$ here $coA$ denotes
a convex hull of the set $A$.
\end{thm}

\begin{pf}  The first part of the proof is similar to the proof of
the theorem \ref{4}.

Therefore, consider when  $l=2$.  From Lemma \ref{3} it follows
that $2p_{ij,1}\geq 1$ and $2p_{i2,2}\geq 1$, hence
$p_{i2,1}+p_{i2,2}\geq 1$. But $\sum\limits_{k=1}^mp_{i2,k}=1$
implies  that $2p_{i2,1}=2p_{i2,2}=1$ and $p_{i2,k}=0$ for all
$k\geq 3.$ Now we can rewrite operator (5) as:
$$
\left.
\begin{array}{lll}
(Vx)_1=x_1+\sum\limits_{i=3}^mx_i^2+2\sum\limits_{1<i<j}p_{ij,1}x_ix_j\\[5mm]
(Vx)_2=x_2+2\sum\limits_{1<i<j\ \ i\neq 2}p_{ij,2}x_ix_j\\[5mm]
(Vx)_k=2\sum\limits_{1<i<j}p_{ij,k}x_ix_j, \  \ 3\leq k\leq m
\end{array}
\right\}\eqno(8)
$$

Putting $x_{\lambda}=\lambda e_1+(1-\lambda)e_2 $, where $0\leq
\lambda\leq 1$, we get $Vx_{\lambda}=x_{\lambda}.$ Therefore, $V$
has infinitely many fixed points $x_\lambda$. Simple calculations
show that Jacobian $J(x_{\lambda})$ of the fixed point
$x_{\lambda}$ has a following view:

$$
\left(%
\begin{array}{ccccc}
  1 & 0 & 0 & ... & 0 \\
  0 & 1 & 0 & ... & 0 \\
  0 & 0 & 0 & ... & 0 \\
  ... & ... & ... & ... & ... \\
  0 & 0 & 0 & ... & 0 \\
\end{array}%
\right)$$

It is easy to see that the only eigenvalue of this matrix is 1,
which in its turn belong to the unit ball. Therefore, all of fixed
points are hyperbolic.

 Now consider a function
$\varphi:S^{m-1}\rightarrow R,$ defined by
$\varphi(x)=x_3+x_4+\cdots+x_m.$ Then for any $x\in S^{m-1}$ we
have
$$\varphi(Vx)=\sum\limits_{1<i<j}\sum\limits_{k=3}^{m}2p_{ij,k}x_ix_j.$$

The inequality $2p_{ij,1}-1\geq 0$ implies
$2\sum\limits_{k=2}^mp_{ij,k}\leq 1$, which yields
$$\varphi(Vx)\leq \sum\limits_{1<i<j}x_ix_j\leq
\sum\limits_{i=3}^mx_i\sum\limits_{i=1}^mx_i=\sum\limits_{i=3}^mx_i=\varphi(x).$$
Consequently, $\{\varphi(V^nx)\}$ is a decreasing sequence.
Therefore it converges.

Denote
$$\lim\limits_{n\rightarrow \infty}\varphi(V^nx)=C.
$$

From (8) one gets
\begin{eqnarray*}
(V^{n+1}x)_1+(V^{n+1}x)_2&=&(V^{n}x)_1+(V^{n}x)_2+2\sum\limits_{1<i<j}p_{ij,1}(V^{n}x)_i(V^{n}x)_j\\
&&+2\sum\limits_{1<i<j}p_{ij,2}(V^{n}x)_i(V^{n}x)_j.
\end{eqnarray*}

According to
$$\lim\limits_{n\rightarrow
\infty}((Vx)_1^{(n+1)}+(Vx)_2^{(n+1)})=\lim\limits_{n\rightarrow
 \infty}((Vx)_1^{(n)}+(Vx)_2^{(n)})=1-C$$
we have
 $$\lim\limits_{n\rightarrow\infty}
 (\sum\limits_{1<i<j}p_{ij,1}(Vx)_i^{(n)}(Vx)_j^{(n)}+
\sum\limits_{1<i<j}p_{ij,2}(Vx)_i^{(n)}(Vx)_j^{(n)})= 0,$$ which
means that $C=0$, and therefore $\omega (x)\in co\{e_1,e_2\}.$

\textbf{Remark.} We showed that all of the fixed points are
parabolic, this means that trajectory in a neighborhood of these
points is nonstable.

\textbf{Observation} In last two theorems we saw that the
trajectory of initial point tend to the bound of the simplex.
therefore it is natural to ask whether dissipative q.s.o. and
Volterra q.s.o coincide. An answer is negative

Dissipative q.s.o. (5) can be Volterra q.s.o. if and only if
$\alpha_k=\{k\}.$ On the other hand  from the lemma 3 it
automatically follows that $p_{ik,k}=\frac12$ and
$a_{ki}=2p_{ik,k}-1=0$ Therefore only identity operator can be
contemporary dissipative and Volterra q.s.o. .
\end{pf}

\section{Ergodicity}

In this section we are going to show that any dissipative q.s.o.
is ergodic.

Recall that a q.s.o. is called {\it ergodic} if the following
limit exists
$$
\lim\limits_{n\rightarrow\infty}\frac{x+Vx+\cdots+V^{n-1}x}{n}
$$
for any $x\in S^{m-1}$.

Ulam \cite{Ul} formulated a conjecture that \textit{any q.s.o. is
ergodic.} However, Zakharevich \cite{Z} showed that it is not so.
He considered the following q.s.o.
\begin{eqnarray*}
&&(Vx)_1=x_1^2+2x_1x_2,\\
&&(Vx)_2=x_2^2+2x_2x_3,\\
&&(Vx)_3=x_3^2+2x_1x_3
\end{eqnarray*} and proved that such an operator is not
ergodic.

Now we show that Ulam's conjecture is true for dissipative q.s.o.

\begin{thm}\label{6} Any dissipative q.s.o. is
ergodic.
\end{thm}

\begin{pf} Let $V:S^{m-1}\rightarrow S^{m-1}$ be a dissipative
q.s.o. Then we have
$$x\prec Vx\prec V^2x\prec V^3x\prec \cdots $$
It means that
\begin{eqnarray*}
&&x_{[1]}\leq (Vx)_{[1]}\leq (V^2x)_{[1]}\leq\cdots\\
&&x_{[1]}+x_{[2]}\leq (Vx)_{[1]}+(Vx)_{[2]}\leq
(V^2x)_{[1]}+(V^2x)_{[2]}\leq\cdots\\
&&\vdots\\
&&\sum\limits_{i=1}^kx_{[i]}\leq \sum\limits_{i=1}^k(Vx)_{[i]}\leq
\sum\limits_{i=1}^k(V^2x)_{[i]}\leq\cdots
\end{eqnarray*}

The sequences $\{\sum\limits_{i=1}^k(V^n(x))_{[i]}, \ \
n=1,2,\cdots\} \ \forall k=\overline{1,m}$ are increasing and
bounded, consequently convergent. The last means that the
following sequences are also convergent
$$\{(V^n(x))_{[k]} \ \ n=1,2,\cdots \} \ \forall k=\overline {1,m}.$$

Let's denote $y_k=\lim\limits_{n\rightarrow \infty}(V^nx)_{[k]},$
and $y=(y_1,y_2,\cdots,y_m)$

If $z=(z_1,z_2,\cdots,z_m)\in \omega(x^0),$ then there exists
$\{x^{(n_j)}\},$ such that $(V^{n_j}x)\rightarrow z.$ Therefore we
have $(V^{n_j}x)_{\downarrow}\rightarrow z_{\downarrow}.$ On the
other hand $(V^{n_j}x)_{\downarrow}\rightarrow y,$ since $y$ is a
limit of the sequence $(V^nx), \ n=1,2,\cdots$  That's why
$z_{\downarrow}=(y_1,y_2,\cdots,y_m)=y.$ Therefore, we infer that
any element of $\omega(x)$ is some kind of rearrangement of
$(y_1,y_2,\cdots,y_m)=y$. This means that the
 cardinality of $\omega(x)$ cannot be greater than $m!$

Let $|\omega(x)|=p,$ then the trajectory of the $\{V^nx\}$ tends
to the cycle of order $p$, i.e. the trajectory is divided into $p$
convergent subsequences. The operator $V$ acts as a cyclic
permutation of their limits. Therefore we conclude that $V$ is
ergodic.
\end{pf}

\textbf{Remark.} Now we proved that limit set of the trajectory is
finite. From biological point of view this means that there are
periodical evolutions, since there are periodical points of the
dissipative q.s.o.

\textbf{Remark.} Let $f:K\rightarrow K,$ here $f$ is continuous
and $K$ Hausdorff space. It is well known fact that if $\omega$
limit set of trajectory of any initial point is finite, then it is
ergodic and for any subsequence $n_k\in \texttt{N}$ is also
ergodic. That is the following limit exists
$$
\lim\limits_{k\rightarrow\infty}\frac{V^{n_1}x+V^{n_2}x+\cdots+V^{n_k}x}{k}
$$

\section{Dissipative Cubic Stochastic Operators.}

Recently in \cite{R} a notion of cubic stochastic operator was
introduced and studied a class of such kind of operators. Namely,
In the similar manner, the following results can be proved.

Let us denote $a_{ijl}=(p_{ijl,1},p_{ijl,2},\cdots, p_{ijl,m})\ \
\forall i,j,l=\overline {1,m}$ where $p_{ijl,k}$ are coefficients
of c.s.o. (3).

\begin{lem}\label{7} If $V$ is a dissipative c.s.o. then
$$ (a_{iii})_{\downarrow} =e_1 \ \ \forall i=\overline {1,m}.$$
\end{lem}

From Lemma \ref{7} we can rewrite dissipative c.s.o. in the
following form
$$
(Vx)_k=\sum\limits_{i\in
\alpha_k}x_i^2+\sum\limits_{(i-j)^2+(j-l)^2+(l-i)^2>0}p_{ijl,k}x_ix_jx_l.\eqno(9)
$$
Here $\alpha_k\subset I=\{1,2\cdots,m\}, \ \ \ \alpha_i\cap
\alpha_j=\emptyset, \ \ i\neq j \ \ \
\bigcup\limits_{k=1}^m\alpha_k=I.$

\begin{lem}\label{8} Let (9) be a dissipative c.s.o. If $j\in
\alpha_{k_0},$ then $p_{ijj,k_0}=(a_{ijj})_{[1]}\geq \frac{2}{3} \
\ \forall i=\overline{1,m}$ and $(a_{ijj})_{[k]}=0 \ \ \forall
k\geq 3, \ \forall i=\overline{1,m}$
\end{lem}

Note that necessary conditions for dissipativity of c.s.o. is not
solved completely.

\begin{thm} Any dissipative c.s.o. is  ergodic. \end{thm}

\textbf{Remark.} In proved ergodicity we used only the condition
(4) and finite dimensionality of the space, therefore we can
deduce that any dissipative stochastic operator is ergodic.

\section{Conclusion.} The main achievement of the present paper
is creating and studying a new class of the quadratic stochastic
operators. Dissipative q.s.o. has various application in
mathematical genetics and one can use given results. The main
results of this work are theorems 3.3, 3.5, and 4.1 . The methods,
which was used for proving results are different from those
well-known methods. One can use these methods and techniques for
proving another results, not only, in the theory of quadratic
stochastic operators, but also in other disciplines of
mathematics, namely nonlinear analysis, dynamical systems and
ergodic theory. Nevertheless, in class of dissipative q.s.o. there
are some open problems.

Prove or disprove the following statements.

\textbf{Problem 1.} Any dissipative q.s.o. has either unique or
infinitely many fixed points.

\textbf{Problem 2.} If $V$ is a dissipative q.s.o., then $\omega-$
limit set of any non-fixed initial belongs to the bound of the
simplex.

\textbf{Problem 3.} Now recall that c.s.o. (3) is said to be
Volterra c.s.o. if $p_{ijl,k}=0 \ \forall k\notin \{i,j,l\}.$
Whether there is non identity dissipative Volterra c.s.o.?

\section*{Acknowledgements}

I am very thankful to Prof. R.N. Ganikhodzhaev and Prof. F.M.
Mukhamedov for fruitful discussions and encourage to present work.
The work also partially supported by ICTP, OEA-AC-84.

\end{document}